\documentclass[12pt]{article}
\usepackage{amsfonts,amsmath,amsthm,amssymb}
\usepackage{graphicx}
\usepackage{cite}
\title{A Geometric(1/2) Distribution Associated with Record Breaking}
\author{
        Dan Naiman \& Fred Torcaso \\
                Department of Applied Mathematics and Statistics\\
                Johns Hopkins University\\
                Baltimore, MD
 }
\date{\today}

\begin{document}

\maketitle
\begin{abstract}
Let $X_i,i=0,1,\ldots$ be a sequence of iid random variables whose distribution is continuous.
Associated with this sequence is the sequence $(i,X_i),i=0,1,\ldots$.
Let ${\cal R}_{n}$ denote the set of Pareto optimal elements of $\{ (i,X_i):i=0,\ldots,n\}.$
We refer to the elements of ${\cal R}_{n}$ as the \emph{current records at time} $n,$ and we define $R_n=\vert {\cal R}_n\vert,$
the number of such records.  Observe that $R_n$ has $\{1,\ldots,n+1\}$ as its support.
When $(n,X_n)$ is realized, it is a Pareto optimal element of $\{ (i,X_i)~:~i=0,\ldots,n\}$ and
${\cal R}_{n} \backslash (n,X_n) \subset {\cal R}_{n-1}.$
Then we refer to those elements of ${\cal B}_n = {\cal R}_{n-1} \backslash {\cal R}_{n}$ as the records \emph{broken} at time $n.$
Let $B_n= \vert {\cal B}_n \vert.$ We show that
$$
P[B_n = k] \rightarrow 1/2^{k+1} \mbox{ for } k=0,1,2,\ldots.
$$
\end{abstract}

\section{Introduction}
Let $X_i,i=0,1,\ldots$ be a sequence of iid random variables whose distribution is continuous.
Associated with this sequence is the sequence $(i,X_i),i=0,1,\ldots$.
We introduce the partial ordering on such pairs by $(i,X_i) \prec (j,X_j)$ if $i\leq j$ and $X_i \leq X_j.$
Let ${\cal R}_{n}$ denote the set of Pareto optimal elements of $\{ (i,X_i):i=0,\ldots,n\}.$
These are the pairs $(i,X_i)$ for $i=1,\ldots,n$ for which $(i,X_i) \not\prec (j,X_j)$ for all $1 \leq j\leq n$ with $j \neq i.$

\begin{figure}
\begin{center}
\includegraphics[scale=.75]{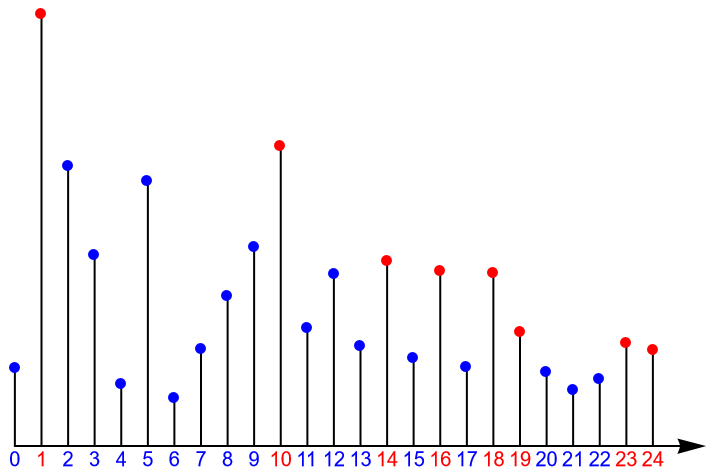}
\end{center}
\caption{A plot of $(i,X_i),~$ $i=0,\ldots,24$ for a  particular realization.\\ The red points are the points on the Pareto boundary of $\{(i,X_i),i=0,\ldots,24\}$ so
 ${\cal R}_{24}=\{(1,X_1),(10,X_{10}),(14,X_{14}), (16,X_{16}), (18,X_{18}), (19,X_{19}),(23,X_{23}), (24,X_{24})\}$
and $R_{24}=8.$ The current records at time 24 correspond to the records in the classical sense for the reverse sequence $X_{24},X_{23},\ldots,X_0$}.
\end{figure}

We refer to the elements of ${\cal R}_{n}$ as the \emph{current records at time} $n,$ and we define $R_n=\vert {\cal R}_n\vert.$
Observe that at time $n$ the \emph{current records} correspond to those indices $i$ in which univariate records, in the traditional sense, occur for the reversed sequence $X_{n},\ldots,X_{0}.$  (See Figure 1).
the number of such records. (See Figure 1.) Observe that $R_n$ has $\{1,\ldots,n+1\}$ as its support.
When $(n,X_n)$ is realized, it is a Pareto optimal element of $\{ (i,X_i)~:~i=0,\ldots,n\}$ and
${\cal R}_{n} \backslash (n,X_n) \subset {\cal R}_{n-1}.$
Then we refer to those elements of ${\cal B}_n = {\cal R}_{n-1} \backslash {\cal R}_{n}$ as the records \emph{broken} at time $n.$
Let $B_n= \vert {\cal B}_n \vert.$  (See Figure 2.)

\begin{figure}
\begin{center}
\includegraphics[scale=.75]{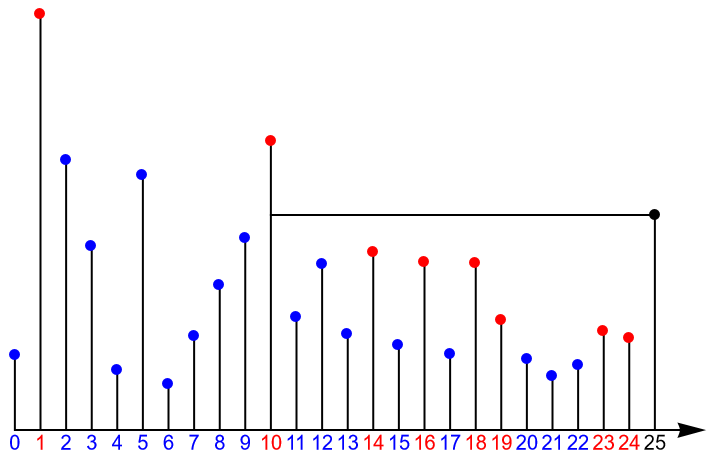}
\end{center}
\caption{When $n=25$ the new set of current records is ${\cal R}_{25}=\{(1,X_1),(10,X_{10}),(25,X_{25}\}$ and $B_{25}=6$ records are broken at this time
}
\end{figure}

Our main result is the following.

{\bf Theorem.}
The limiting distribution of $B_n$ is Geometric(1/2), that is,
$$
P[B_n = k] \rightarrow 1/2^{k+1}   \mbox{ for } k=0,1,2,\ldots.
$$

For the situation when $(X_i,Y_i)$ are iid with $X_i$ and $Y_i$ independent with continuous distributions, one of the authors observed empirically that, conditioned on $(X_n,Y_n)$ being a record in the Pareto optimal sense, the number of current records it breaks has the same Geometric(1/2) limiting distribution, and this observation has been confirmed as a theorem by Fill \cite{fill}. While the conclusions of these two results are similar, it is not clear why these two models give rise to similar results.

{\bf Some Preliminaries.}

We'll assume, without loss of generality that $X_i \sim \mbox{Uniform}(0,1)$ and since the probability of any tie among the $X_i$ is zero, so we will assume that ties do not occur.
We let $C_{n-1}$ denote the set of indices for the current records at time $n-1,$ so that $R_{n-1} = \vert C_{n-1} \vert.$
Assuming $R_{n-1} = r,$ let the indices in $C_{n-1}$ in reverse order be denoted by $I_0^{(n-1)},I_1^{(n-1)},\ldots,I_{r-1}^{(n-1)}.$
As an immediate consequence we see that

\begin{itemize}
\item $I_{r-1}^{(n-1)} < I_{r-2}^{(n-1)} < \cdots < I_1^{(n-1)} < I_0^{(n-1)}=n-1,$
\item $X_{I_{r-1}} > X_{I_{r-2}} > \cdots > X_{I_1} > X_{I_0},$ and
\item $X_j < X_{I_p}$ for $I_{p+1}<j < I_p$
\end{itemize}

In the example shown in Figure 1, we have
\begin{itemize}
\item $C_{24} = \{ 1,10,14,16,18,19,23,24\}.$
\item $r=8$
\item $(I_{r-1}^{(24)},I_{r-2}^{(24)},\cdots, I_1^{(24)},I_0^{(24)})=(1,10,14,16,18,19,23,24).$
\end{itemize}

{\bf Proof of Theorem.}

For $k=0$ the event $\{ B_n = 0 \}$ is equivalent to $\{ X_{n}<X_{n-1}\}$ which has probability 1/2.

For $\{B_n=1\}$ to occur, it must be the case that there is at least one current record at time $n-1$ and we have
$$
P[B_n=1]= P[ B_n=1, R_{n-1}=1] + P[ B_n=1, R_{n-1}\geq 2].
$$

The event $\{B_n=1, R_{n-1}=1\}$ is equivalent to the order statistics satisfying $X_{(n-1)}=X_{n-1}$ and $X_{(n)}=X_n$ and the probability of this event is
$\frac{1}{n(n+1)}.$

For the second term, note that when exactly one record is broken at time $n$ the broken record is at index $I_0^{(n-1)}=n-1$ and the last unbroken record is
at $I_1^{(n-1)}.$
We decompose the second event according to the value of $I^{(n-1)}_1$  to give
$$
P[B_n=1]= P[ B_n=1, R_{n-1}=1] + \sum_{i=0}^{n-2} P[ B_n=1, R_{n-1}\geq 2, I_1^{(n-1)}=i].
$$

For any fixed $i=0,\ldots,n-2$ we have (see Figure 3)
$$
P[B_n=1, R_{n-1}\geq 2, I_1^{(n-1)}=i]=
P[\bigcap_{j=i+1}^{n-2}  \{X_j <X_{n-1}\}\cap \{X_{n-1} < X_n<X_i\}]
$$

\begin{figure}
\begin{center}
\includegraphics[scale=.75]{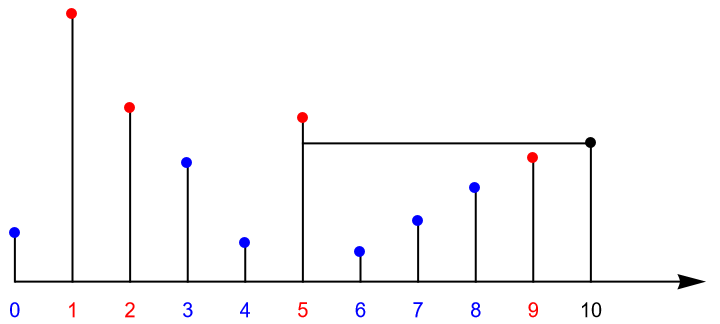}
\end{center}
\caption{When a single record is broken at time $n$ that record must be $(n-1,X_{n-1})$ as this is always a current record at time $n-1$ so $I_0^{(n-1)}=n-1.$
In the case shown, a single record is broken at time $10,$ and the index of the last unbroken record is $I^{(n-1)}_1=5.$ Hence $X_j<X_9.$ for $j=6,7,8$ and $X_9<X_{10}<X_5.$
}
\end{figure}

and conditioning on $(X_i,X_{n-1})=(v,u)$ this becomes
$$
= \int_{0}^1 dv \int_{0}^v u^{n-i-2}(v-u) \,du
$$
$$
= \int_{0}^1 dv \int_{0}^v u^{n-i-2}v\,du -\int_{0}^1 dv \int_{0}^v u^{n-i-1}\,du
$$
$$
= \int_{0}^1 v^{n-i}/(n-i-1)\,dv -\int_{0}^1 v^{n-i}/(n-i)\,dv
$$
$$
=\left( \frac{1}{n-i-1}-\frac{1}{n-i}\right) \int_{0}^1 v^{n-i}\,dv
$$
$$
=\frac{1}{(n-i-1)(n-i)(n-i+1)}
$$
Since
$$
\sum_{i=1}^n \frac{1}{i(i+1)(i+2)} = \frac{1}{4} - \frac{1}{2(n+1)(n+2)}
$$
we find that
$$
P[B_n=1] = \frac{1}{n(n+1)} + \sum_{i=0}^{n-2} \frac{1}{(n-i-1)(n-i)(n-i+1)}
$$
$$=
\frac{1}{n(n+1)} +\sum_{j=1}^{n-1} \frac{1}{j(j+1)(j+2)}
$$
$$
=\frac{1}{4} + \frac{2}{n(n+1)}
$$

For the general case, if $k>2,$ we have
$$
P[B_n=k] = P[B_n=k, R_{n-1}\leq k]+P[B_n=k,R_{n-1}\geq k+1].
$$
Since $R_n$ is the number of records in the classical sense for the reverse sequence $X_{n-1},\ldots,X_0$
we have $R_n/\log(n)\stackrel{a.s.}{\rightarrow} 0$
(see \cite{Arnold} p. 25)
making the first probability $o(1).$

For the second probability, under the assumtion that $R_{n-1}\geq k+1$ the indices $I^{(n-1)}_{k},\ldots,I^{(n-1)}_{0}=n-1$ are well-defined and form increasing sequence,
and the event $\{ B_n = k,R_{n-1}\geq k+1\}$ is equivalent to the existence of indices $0\leq i_k<i_{k-1}<\cdots<i_1 < i_0=n-1,$ with the properties that
\begin{itemize}
\item[(a)] $I_p=i_p,$ for $p=0,\ldots,k,$
\item[(b)] $X_{i_k}> X_{i_{k-1}}> \cdots > X_{i_0},$
\item[(c)] $X_j < X_{i_{p-1}}$ for $j=i_p+1,\ldots,i_{p-1}-1,$ for $p=1,\ldots,k,$ and
\item[(d)] $X_{i_{k-1}} < X_n < X_{i_k}.$
\end{itemize}

Let
$${\cal I}_{j} = \{(i_1,\ldots,i_j) : k-p \leq i_p \leq i_{p-1}-1 \mbox{ for } p=1,\ldots,j\}.$$
So
\begin{equation}
\label{probability0}
P[B_n=k,R_{n-1}\geq k+1] =
\sum_{(i_1,\ldots,i_k)\in {\cal I}_k}
p(i_0,\ldots,i_k),
\end{equation}
where $p(i_0,\ldots,i_k)$ denotes the probability of the event
$$
\left(\bigcap_{p=1}^{k} \{ X_{i_p}>X_{i_{p-1}}\}\right)
 \cap
\left( \bigcap_{p=1}^{k-1}
 \bigcap_{j=i_{p-1}+1}^{i_p+1} \{X_j < X_{i_{p-1}-1}\}\right)
 \cap \{X_{i_{k-1}} < X_n < X_{i_k}\}.
$$
Conditioning on $X_{i_0}=u_0,\ldots,X_{i_k}=u_k$ we see that

\noindent $p(i_0,\ldots,i_k)$

\noindent \hskip .1in $\displaystyle
= \int_{0}^1\!\! du_k\!\!\int_{0}^{u_k}\!\!du_{k-1} \cdots \!\! \int_{0}^{u_2} \!\!du_1\!\!\int_{0}^{u_1}
\!\! (u_k-u_{k-1})\prod_{p=0}^{k-1} u_p^{i_p-i_{p+1}-1}\,du_0
$\vskip -.1 in

\noindent \hskip .1 in $\displaystyle
= \frac{1}{i_0-i_1} \int_{0}^1\!\! du_k\!\!\int_{0}^{u_k}\!\!du_{k-1} \cdots  \!\! \int_{0}^{u_2}
\!\! (u_k-u_{k-1})\prod_{p=2}^{k-1} u_p^{i_p-i_{p+1}-1} u_1^{i_0-i_2-1}\,du_1
$\vskip -.1 in

\noindent \hskip .1in $\displaystyle
=  \frac{1}{i_0-i_1}\frac{1}{i_0-i_2} \int_{0}^1 \!\!du_k\!\!\int_{0}^{u_k}\!\!du_{k-1}  \cdots \!\! \int_{0}^{u_3}
\!\! (u_k-u_{k-1})\prod_{p=3}^{k-1} u_p^{i_p-i_{p+1}-1} u_2^{i_0-i_3-1}\,du_2
$

\noindent \hskip .1in $
\vdots
$

\noindent \hskip .1in $\displaystyle
=  \left(\prod_{p=1}^{k-2} \frac{1}{i_0-i_p}\right)
\int_{0}^1 \!\!du_k\!\!\int_{0}^{u_k}\!\!du_{k-1} \!\! \int_{0}^{u_{k-1}}
\!\! (u_k-u_{k-1}) u_{k-1}^{i_{k-1}-i_{k}-1}
 u_{k-2}^{i_0-i_{k-1}-1}\,du_{k-2}
$\vskip -.1in

\noindent \hskip .1in $\displaystyle
=  \left(\prod_{p=1}^{k-1} \frac{1}{i_0-i_p}\right)
\int_{0}^1 \!\!du_k\!\! \int_{0}^{u_k}
 (u_k-u_{k-1}) u_{k-1}^{i_0-i_k-1}\,du_{k-1}
$\vskip -.1 in

\noindent \hskip .1in $\displaystyle
=\left(  \prod_{p=1}^{k-1} \frac{1}{i_0-i_p}\right) \left( \frac{1}{(i_0-i_k)(i_0-i_k+2)}-\frac{1}{(i_0-i_k+1)(i_0-i_k+2)}\right)
$\vskip -.1 in

\noindent \hskip .1in $\displaystyle
=\left(  \prod_{p=1}^{k-1} \frac{1}{i_0-i_p}\right)
\left( \frac{1}{(i_0-i_k)(i_0-i_k+1)(i_0-i_k+2)}\right).
$

\noindent Substituting this expression in (\ref{probability0}) we obtain
\begin{equation}
\label{prob}
P[B_n=k, R_{n-1}\geq k+1]
\end{equation}
\begin{equation}
\label{probabilityexpression}
= \sum_{(i_1,\ldots,i_{k})\in I_{k}}
\left(  \prod_{p=1}^{k} \frac{1}{i_0-i_p}\right)
 \frac{1}{(i_0-i_k+1)(i_0-i_k+2)}
\end{equation}
$$
=\sum_{(i_1,\ldots,i_{k-1})\in I_{k-1}}
\left(  \prod_{p=1}^{k-1} \frac{1}{i_0-i_p}\right)
\sum_{i_k=0}^{i_{k-1}-1} \frac{1}{(i_0-i_k)(i_0-i_k+1)(i_0-i_k+2)}.
$$
\begin{equation}
\label{probability1}
=
\sum_{(i_1,\ldots,i_{k-1})\in I_{k-1}}
\left(  \prod_{p=1}^{k-1} \frac{1}{i_0-i_p}\right)
\sum_{j=i_0-i_{k-1}+1}^{i_0} \frac{1}{j(j+1)(j+2)}
\end{equation}

Writing
$$
\sum_{j=i_0-i_{k-1}+1}^{i_0} \frac{1}{j(j+1)(j+2)} = \sum_{j=1}^{i_0} \frac{1}{j(j+1)(j+2)}
-\sum_{j=1}^{i_0-i_{k-1}} \frac{1}{j(j+1)(j+2)}
$$
and using the identity
$$
\sum_{i=1}^m \frac{1}{i(i+1)(i+2)} = \frac{1}{4} - \frac{1}{2(m+1)(m+2)}
$$
(\ref{probability1}) takes the form
\begin{equation}
\label{sumexpression1}
\sum_{(i_1,\ldots,i_{k-1})\in I_{k-1}}
\left(  \prod_{p=1}^{k-1} \frac{1}{i_0-i_p}\right)
\left\{  \frac{1}{2(i_0-i_{k-1}+1)(i_0-i_{k-1}+2)} \right\}
\end{equation}

\begin{equation}
\label{sumexpression2}
- \sum_{(i_1,\ldots,i_{k-1})\in I_{k-1}}
\left(  \prod_{p=1}^{k-1} \frac{1}{i_0-i_p}\right)
\left\{ \frac{1}{2(i_0+1)(i_0+2)} \right\}.
\end{equation}
The second sum in expression (\ref{sumexpression2}) can be written as
\begin{equation}
\label{secondterm}
\frac{1}{2(i_0+1)(i_0+2)}\sum_{i_1=k-1}^{i_0-1}\frac{1}{i_0-i_1}
\sum_{{i_2}=k-2}^{i_1-1} \frac{1}{i_0-i_2}
 \cdots
\sum_{i_{k-1}=1}^{i_{k-2}-1} \frac{1}{i_0-i_{k-1}}.
\end{equation}
and each sum starting with the inner-most one can be bounded above. Indeed, for any of the sums we have
$$
\sum_{i_p=k-p}^{i_{p-1}-1} \frac{1}{i_0-i_p} = \sum_{j=i_0-i_{p-1}+1}^{i_0-k+p}\frac{1}{j}
= \sum_{j=1}^{i_0-k+p}\frac{1}{j} - \sum_{j=1}^{i_0-i_p}\frac{1}{j} \leq 1+\log(i_0-k+p),
$$
for $p=1,\ldots,k-1.$
Recalling that $i_0=n-1$ and the term (\ref{secondterm}) (which is nonnegative) is bounded above by
$$
\frac{(1+\log(n-1))^{k-1}}{2n(n+1)} = o(1).
$$
We conclude that  (\ref{prob}) equals
$$
\frac{1}{2}\sum_{(i_1,\cdots,i_{k-1}) \in {\cal I}_{k-1}}
\left(  \prod_{p=1}^{k-1} \frac{1}{i_0-i_p}\right)
\left\{  \frac{1}{(i_0-i_{k-1}+1)(i_0-i_{k-1}+2)} \right\} + o(1).
$$

The only differences between this expression  and (\ref{probabilityexpression}) are
\begin{itemize}
\item[(a)] the index set being summed over changes from ${\cal I}_k$ to ${\cal I}_{k-1},$
\item[(b)] a factor of $\frac{1}{2}$ is introduced, and
\item[(a)] an additional $o(1)$ term appears.
\end{itemize}

Repeating this reduction process an additional $k-2$ times leads to the conclusion that
 (\ref{prob}) equals
$$
\frac{1}{2^{k-1}} \sum_{i_1=k-1}^{i_0-1}
\frac{1}{(i_0-i_1)(i_0-i_1+1)(i_0-i_1+2)} + o(1)
$$
$$
= \frac{1}{2^{k-1}} \sum_{j=1}^{i_0-k+1} \frac{1}{j(j+1)(j+2)} + o(1)
$$
$$
=  \frac{1}{2^{k-1}} \frac{1}{4} + \frac{1}{2^{k-1}} \frac{1}{2(i_0-k+2)(i_0-k+3)} + o(1)
$$
$$
= \frac{1}{2^{k+1}} + o(1).
~~~~~\Box$$

\bibliography{references}{}
\bibliographystyle{plain}
\end{document}